\documentclass[dvips,preprint,12pt]{imsart}
\RequirePackage[OT1]{fontenc}
\RequirePackage{amsthm,amsmath,natbib}
\RequirePackage[colorlinks]{hyperref}
\RequirePackage{hypernat}
\usepackage{amsmath,epsfig}
\usepackage{amstext,amssymb,amsmath,amscd,color}
\usepackage{dsfont}
\startlocaldefs
\def\e{\rm e}
\def\bea{\begin{eqnarray}}
\def\eea{\end{eqnarray}}
\def\beas{\begin{eqnarray*}}
\def\eeas{\end{eqnarray*}}
\numberwithin{equation}{section}
\theoremstyle{plain}
\newtheorem{theorem}{Theorem}[section]
\newtheorem{definition}{Definition}[section]

\endlocaldefs
\def\be{\begin{equation}}
\def\ee{\end{equation}}
\begin{document}
\begin{frontmatter}
\title{Smooth Kernel Estimation of a Circular Density Function: A Connection to Orthogonal Polynomials on the Unit Circle}
\runtitle{Smooth Circular Density: A Connection to OPUC}
\runauthor{Y.P. Chaubey}
\begin{aug}
\author{\fnms{} \snm{Yogendra P. Chaubey}\thanksref{}\ead[label=c]
{yogen.chaubey@concordia.ca}}
\address{Department of Mathematics and Statistics, Concordia University,\\
Montr\'eal, QC H3G 1M8 \\
\printead{c}}
\end{aug}

\begin{abstract}
In this note we provide a simple approximation theory motivation for the circular kernel density estimation and further explore the usefulness of the wrapped Cauchy kernel in this context. It is seen that the wrapped Cauchy kernel appears as a natural candidate in connection to orthogonal series density estimation on a unit circle. This adds further weight to the considerable role of the wrapped Cauchy in circular statistics.
\end{abstract}


%\tableofcontents
\end{frontmatter}
\thispagestyle{empty}
\section{Introduction}
Consider an absolutely continuous (with respect to the Lebesgue measure) circular density $f(\theta),\theta \in [-\pi,\pi],$ {\it i.e} $f(\theta)$ is $2\pi-$periodic,
\be\label{circden}
f(\theta)\ge 0~{\rm for}~\theta\in \mathds{R} ~{\rm and}~\int_{-\pi}^{\pi}f(\theta)d\theta=1.
\ee
In the literature on modeling circular data, starting from the classical text of \cite{mardia72}, there appear many standard texts such as \cite{fisher93}, \cite{jammalSen01} and \cite{mardiajupp00} that cover parametric models along with many inference problems. More recently  various alternatives to these classical parametric models, exhibiting asymmetry and multimodality  have been investigated with respect to their mathematical properties and goodness of fit to some real data; see \cite{AbePewsey11}, \cite{jonespewsey12}), \cite{katojones15}, \cite{katojones10}, \cite{minhfarnum03} and \cite{Shimizuida02}.

In cases where multimodal and\/or asymmetric models may be appropriate, semiparametricr or nonparametric modelling may be considered more appropriate. \cite{fernduran04} and \cite{mooneyetal03} considered semi-parametric analysis based on mixture of circular normal and von Mises distributions and \cite{halleta87}, \cite{Baietal88}, \cite{fisher89}, \cite{Taylor08} and \cite{klemela00}  have considered nonparametric approaches.

Given a random sample $(\theta_1,....\theta_N)$ from the density \eqref{circden}, the circular kernel density estimator is given by
\be\label{npkerden}
\hat f(\theta;h)=\frac{1}{N}\sum_{j=1}^{N}k_h(\theta-\theta_j)
\ee
where $k_h(\theta-\phi)$ is a circular kernel density function that is concentrated around $\phi$ as $h\to h_0$ for some known $h_0.$ As motivated in \cite{Taylor08}, a natural choice for the kernel function is one of the commonly used circular probability densities, such as the wrapped normal distribution, or the von Mises distribution. \cite{Taylor08} investigated the use of von Mises kernel, in which case the density estimator is given by
\be\label{vMkernelden}
\hat f_{VM}(\theta;\nu)=\frac{1}{N(2\pi)I_0(\nu) }\sum_{j=1}^{N}\exp\{\nu(\theta-\theta_j)\}
\ee
where $I_0(\nu)$ is the Bessel function of order $r$ and $\nu$ is the concentration parameter. \cite{dimaretal09} considered the use of circular kernels to circular regression while extending the use of von Mises kernels to more general circular kernels. In the present note I demonstrate that the wrapped Cauchy kernel presents itself as the kernel of choice by considering an estimation problem on the unit circle. We also show that this approach leads to orthogonal series density estimation, however no truncation of the series is required. It may be noted that the wrapped Cauchy distribution with location parameter $\mu$ and concentration parameter $\rho$ is given by
\be\label{wrappedCauchyden}
f_{WC}(\theta;\mu,\rho)=\frac{1}{{2\pi}}\frac{1-\rho^2}{1+\rho^2-2\rho \cos(\theta-\mu)}~,-\pi \le \theta <\pi,
\ee
that becomes degenerate at $\theta=\mu$ as $\rho \to 1.$ The estimator of $f(\theta)$ based on the above kernel is given by
\be\label{WCkernelden}
\hat f_{WC}(\theta;\rho)=\frac{1}{N}\sum_{j=1}^{N}f_{WC}(\theta;\theta_j,\rho).
\ee

In Section 2, we provide a simple approximation theory argument behind the nonparametric density estimator of the type introduced in \eqref{npkerden} and \eqref{vMkernelden}. In Section 3, first we present some basic results from the literature on orthogonal polynomials on the unit circle and then introduce the strategy of estimating $f(\theta)$ by estimating an expectation of a specific complex function, that in turn produces the non-parametric circular kernel density estimator in \eqref{WCkernelden}.  The next section shows that the  circular kernel density estimator is equivalent to the orthogonal series estimation in a limiting sense. This equivalence establishes a kind of qualitative superiority of the kernel estimator over the orthogonal series estimator that requires the series to be truncated, however the kernel estimator does not have such a restriction.

\section{Motivation for the Circular Kernel Density Estimator}

The starting point of the nonparametric density estimation is the theorem given below from approximation theory (see \cite{mhaskarpai2000}). Before giving the theorem we will need the following definition:
\begin{definition}
Let $\{K_n\}\subset C^*$ where $C^*$ denotes the set of periodic analytic functions with a period $2\pi.$ We say that $\{K_n\}$ is an {\it approximate identity} if \begin{itemize}
\item [A.] $K_n(\theta)\ge 0~ \forall~ \theta\in [-\pi,\pi];$
\item [B.] $\frac{1}{2\pi}\int_{-\pi}^\pi K_n(\theta)=1;$
\item [C.] $\lim_{n\to \infty}\max_{|\theta|\ge \delta}K_n(\theta)=0~{\rm for ~every~}\delta>0.$
 \end{itemize}
\end{definition}

The definition above is motivated from the following theorem which is similar to the one used in the theory of linear kernel estimation (see \cite{prakasarao83}).

\begin{theorem}
Let $f\in C^*,$ $\{K_n\}$ be approximate identity and for $n=1,2,...$ set
\be
f^*(\theta)=\frac{1}{2\pi}\int_{-\pi}^{\pi}f(\eta)K_n(\eta-\theta)d\eta.
\ee Then we have
\be
\lim_{n\to\infty}\sup_{x\in [-\pi,\pi]}|f^*(x) -f(x)|=0.
\ee
\end{theorem}

Note that taking the sequence of concentration coefficients $\rho\equiv \rho_n$ such that $\rho_n\to 1,$ the density function of the Wrapped Cauchy will satisfy the conditions in the definition in place of $2\pi K_n.$ The integral in the above theorem. In general $2\pi K_n$ may be replaced by a sequence of periodic densities on $[-\pi,\pi],$ that converge to to a degenerate distribution at $\theta=0.$

For a  given random sample of $\theta_1,...,\theta_N$ from the circular density $f,$ the Monte-Carlo estimate of $f^*$ is given by
\be
\tilde f(\theta)=\frac{1}{(2\pi)N}\sum_{j=11}^{N} K_n(\theta_j-\theta),
\ee
the suffix $n$ for the kernel $K$ may be a function of the sample size $N.$ The kernel given by the wrapped Cauchy density satisfies the assumptions in the above theorem that provides the estimator proposed in \eqref{WCkernelden}.

\section{Some Preliminary Results from Complex Analysis }

Let $\mathds{D}$ be the open unit disk, $\{z~| ~|z|<1\},$ in $\mathds{Z}$ and let $\mu$ be a continuous measure defined on the boundary $\partial\mathds{D},$ i.e. the circle  $\{z~|~ |z|=1\}.$ The point $z\in \mathds{D}$ will be represented by $z=r{\rm e}^{i\theta}$ for $r\in [0,1), \theta\in [0,2\pi)$ and $i=\sqrt{-1}.$ A standard result in complex analysis involves the Poisson representation that involves the real and complex Poisson kernels that are defined as
\be\label{realPoissonker}
P_r(\theta,\varphi)=\frac{1-r^2}{1+r^2-2r\cos(\theta-\varphi)}
\ee
for $\theta,\varphi\in [0,2\pi)$ and $r\in [0,1)$ and by
\be\label{complexPoissonker}
C(z,\omega)=\frac{\omega+z}{\omega-z}
\ee
for $\omega\in \partial\mathds{D}$ and $z\in \mathds{D}.$ The connection between these kernels is given by the fact that
\be
P_r(\theta,\varphi)={\rm Re} ~C(r\e^{i\theta},\e^{i\varphi})=(2\pi) f_{WC}(\theta;\varphi,\rho).
\ee

The {\it Poisson representation} says that if $g$ is analytic in a neighborhood of $\bar{\mathds{D}}$ with $g(0)$ real, then for $z\in \mathds{D},$

\be\label{Poissonrep}
g(z)=\int \left(\frac{{\e}^{i\theta}+z}{\e^{i\theta}-z} \right){\rm Re}(g(\e^{i\theta}))\frac{d\theta}{2\pi}
\ee
(see \cite[p. 27]{Simon2005}). This representation leads to the result (see (ii) in \S 5 of \cite{Simon2005}) that for Lebesgue a.e. $\theta,$
\be\label{approxF}
\lim_{r\uparrow 1} F(r\e^{i\theta})\equiv F(\e^{i\theta})
\ee
exists and if $d\mu=w(\theta)\frac{d\theta}{2\pi}+d\mu_s$ with $d\mu_s$ singular, then
\be\label{wtheta}
w(\theta)={\rm Re}F(\e^{i\theta}),
\ee
where
\be\label{Fz}
F(z)= \int \left(\frac{{\e}^{i\theta}+z}{{\e}^{i\theta}-z} \right)d\mu(\theta).
\ee
Our strategy for smooth estimation is the fact that for $d\mu_s=0$ we have
\be\label{limitFr}
f(\theta)=\frac{1}{2\pi}\lim_{r\uparrow 1} {\rm Re}~ F(r\e^{i\theta}),
\ee
where

\be
F(z)= \int \left(\frac{\e^{i\theta}+z}{\e^{i\theta}-z} \right)f(\theta)d\theta.
\ee

We define the estimator of $f(\theta)$ motivated by considering an estimator of $F(z),$ the identity \eqref{wtheta} and \eqref{limitFr}, i.e.
\be
\hat f_r( \theta)=\frac{1}{2\pi}{\rm Re}~ F_N(r{\e}^{i\theta})
\ee
where
\be
F_N(r{\e}^{i\theta})=\frac{1}{N}\sum_{j=1}^{N}\left(\frac{{\rm e}^{i\theta_j}+r{\e}^{i\theta}}{{\rm e}^{i\theta_j}-r{\e}^{i\theta}}\right),
\ee
where $r$ has to be chosen appropriately. Recognize that
\be\label{hatF}
F_N(r{\e}^{i\theta})=\frac{1}{N}\sum_{j=1}^{N}~C(z,\omega_j),
\ee
where $\omega_j={\rm e}^{i\theta_j},$ then using \eqref{Poissonrep}, we have
\be
{\rm Re}~F_N(r{\e}^{i\theta})=\frac{1}{N}\sum_{j=1}^{N}~P_r(\theta,\theta_j),
\ee
and therefore
\begin{eqnarray}
\hat f_r( \theta)&=&\frac{1}{(2\pi)N}\sum_{j=1}^{N}~P_r(\theta,\theta_j)\nonumber \\
                 &=&\frac{1}{N}\sum_{j=1}^{N}~f_{WC}(\theta;\theta_j,r),
\end{eqnarray}
that is of the same form as in \eqref{WCkernelden}.
\section{Orthogonal Series Estimation}

We get the orthogonal expansion of $F(z)$ with respect to the basis $\{1,z,z^2,...\}$ as
\be\label{orthoseries}
F(z)=1+2\sum_{n=1}^\infty c_n z^n
\ee
where
$$
c_n=\int {\rm e}^{-in\theta}f(\theta)d\theta,
$$
is the $j^{\rm th}$ trigonometric moment. The series is truncated at some term $N^*$ so that the the error is negligible. However, we show below that estimating the trigonometric moment $c_n,n=1,2,...$ as
$$
\hat c_n=\frac{1}{N}\sum_{j=1}^N{\rm e}^{-in\theta_j},
$$
the estimator of $F(z)$ is the same as given in the previous section. This can be shown by writing
\beas
\hat F(z)&=&1+\frac{2}{N}\sum_{j=1}^N\{ \sum_{n=1}^\infty{\rm e}^{-in\theta_j}z^n\}\\
         &=&1+\frac{2}{N}\sum_{j=1}^N\{ \sum_{n=1}^\infty ({\bar\omega_j}z)^n\};\omega_j={\rm e}^{in\theta_j} \\
         &=&1+\frac{2}{N}\sum_{j=1}^N\left(\frac{{\bar\omega_j}z}{1-{\bar\omega_j}z} \right)\\
         &=&\frac{2}{N}\sum_{j=1}^N\left( \frac{1}{2}+\frac{{\bar\omega_j}z}{1-{\bar\omega_j}z} \right)\\
         &=&\frac{1}{N}\sum_{j=1}^N\left(\frac{1+{\bar\omega_j}z}{1-{\bar\omega_j}z} \right)\\
         &=&\frac{1}{N}\sum_{j=1}^N C(z,\omega_j),
\eeas
which is the same as $F_N(z)$ given in \eqref{hatF}. This ensures that the orthogonal series estimator of the density coincides with the circular kernel estimator.

The determination of the smoothing constant may be handled based on the cross validation method outlined in \cite{Taylor08}.
\bigskip

\noindent{\bf Remark:} Note that the simplification used in the above formulae does not work for $r=1.$ Even though, the limiting form of \eqref{orthoseries} is used to define an orthogonal series estimator
as given by
\be
\hat f_S(\theta)=\frac{1}{2\pi}+\frac{1}{\pi N}\sum_{j=1}^{N}\sum_{n=1}^{n^*} \cos n(\theta-\theta_j),
\ee
where $n^*$ is chosen according to some criterion, for example to minimize the integrated squared error. Thus the above discussion presents two contrasting situations: in one we have to determine the number of terms in the series and in the other number of terms in the series is allowed to be infinite, however, we choose to evaluate ${\rm Re} ~F(r{\rm e}^{i\theta})$ for some $r$ close to 1 as an approximation to ${\rm Re}~ F({\rm e}^{i\theta}).$

\end{document}